\documentclass[letterpaper, 10 pt, conference]{ieeeconf} 
\IEEEoverridecommandlockouts                              
\usepackage{mathtools}
\usepackage{amssymb}
\usepackage{bbold}
\usepackage[yyyymmdd,hhmmss]{datetime}
\usepackage{bm,upgreek}
\usepackage{amssymb}  
\usepackage{color}
\mathtoolsset{showonlyrefs,showmanualtags}
\usepackage{dsfont}
\usepackage{booktabs}
\usepackage[final]{graphicx}
\usepackage{nth}

\DeclareMathOperator{\TR}{Tr}
\DeclareMathOperator{\VAR}{var}

\newcommand{\Exp}[1]{\mathbb{E}\big[ #1\big]}

\usepackage[standard]{ntheorem}
\newtheorem{Assumption}{Assumption}
\newtheorem{Problem}{Problem}

\newcommand{\edit}[1]{\textcolor{black}{#1}}

\title{\LARGE \bf
Near-Optimal Control Strategy in Leader-Follower Networks:  A Case Study for Linear Quadratic Mean-Field Teams 
}

\author{Mohammad M. Baharloo, Jalal Arabneydi and Amir G. Aghdam
\thanks{This work is supported in part by the Natural Sciences and Engineering Research Council of Canada (NSERC) under Grant RGPIN-262127-17, and in part by Concordia University under Horizon Postdoctoral Fellowship.}  
\thanks{Mohammad M. Baharloo, Jalal Arabneydi, and Amir G. Aghdam are with the  Department of Electrical and Computer Engineering, 
        Concordia University, 1455 de Maisonneuve Blvd. West, Montreal, QC, Canada, Postal Code: H3G 1M8.  Email: {\tt\small baharloo@ieee.org}, {\tt\small jalal.arabneydi@mail.mcgill.ca},        
        {\tt\small aghdam@ece.concordia.ca}}%
}

\begin{document}
\maketitle

\vspace*{-5.2cm}{\footnotesize{Proceedings of IEEE Conference on Decision and Control, 2018.}}
\vspace*{4.45cm}

\thispagestyle{empty}
\pagestyle{empty}
\begin{abstract}
In this paper, a decentralized stochastic control system consisting of one leader and many homogeneous followers is studied.  The  leader and followers  are coupled in both dynamics and cost, where the dynamics are linear and the cost function  is quadratic in  the states and actions of the leader and followers.  The objective of the leader and followers is  to  reach consensus while minimizing their communication and energy costs. The leader knows  its local state   and   each follower knows   its local state and  the state of the leader. The number of required links  to implement  this decentralized information structure  is equal to the number of followers, which is the minimum number of links for a  communication graph to be connected. In the special case of leaderless, no link is required among followers, i.e., the communication graph is not even connected.   We  propose a near-optimal control strategy that converges to the optimal solution  as the number of followers increases.  One of the salient features of the proposed solution  is that  it provides a design scheme,  where the  convergence rate  \edit{as well as} the collective behavior of the followers can be designed  by choosing appropriate cost functions.  In addition,  the  computational complexity  of the proposed solution  does not depend on the number of followers.  Furthermore,  the proposed strategy can be computed in a distributed manner, where the leader solves one Riccati  equation and  each follower solves two Riccati equations to calculate their strategies. Two numerical examples are provided to demonstrate the effectiveness of the results in the control of multi-agent systems. 
\end{abstract}

\section{Introduction}
The control of multi-agent systems with leader-follower structure has attracted much interest in the past two decades due to its wide range of applications in various fields of science and engineering. Such applications include vehicle formation~\cite{Fax2004}, sensor networks~\cite{Rawat2014}, surveillance using a team of unmanned aerial vehicles (UAVs)~\cite{UAVsurvey2016} and flocking~\cite{Reynolds1987,Olfati2006flocking}, to name only a few. In this type of problem, a group of agents (called followers) are to track an agent (called leader) while certain performance objectives are achieved. Different performance measures such as minimum energy, fuel or time are considered in the literature. For this purpose, limited communication and computation resources are two main challenges that need to be overcome.

To address the challenges outlined in the previous paragraph, the following two problems have been investigated in the context of consensus control protocols in the literature recently: (i) how the states of the followers can reach the state of the leader under communication constraints (distributed control problem),  \edit{where consensus is reached under  appropriate linear strategies for  properly connected  communication graphs}~\cite{Jadbabaie2003,Ren2007survey,Olfati2007survey}, and (ii) how the state of the followers can reach the state of the leader with minimum energy consumption (optimal control problem),  where   \edit{the optimal control strategy  for  a quadratic performance index  with linear dynamics under centralized information structure is a linear feedback rule obtained by  the solution of the celebrated Riccati equation~\cite{caines1987linear}. Combining the above  two objectives, however, is quite challenging  as it leads to  a decentralized optimal control problem wherein  the optimal} control law is not necessarily linear~\cite{Witsenhausen1968Counterexample}. Furthermore, since the dimension of the matrices in the network model increases with the number of followers, the optimal control law may be intractable for a network of large size. \edit{In this paper, we consider a decentralized optimal control problem with a large number of followers.}

For dynamically decoupled followers and also the case of a leaderless multi-agent system,~\cite{Movric2014,Cao2010} use the control inverse optimality approach to compute the optimal distributed control strategies for special classes of communication graphs. The authors in~\cite{Borelli2008} consider a large number of homogeneous followers and determine the optimal strategy by solving two coupled Riccati equations. 
%
\edit{In contrast, this paper studies a leader-follower multi-agent network with coupled dynamics  under  a directed communication graph in which there is a direct  link from the leader to each follower. In the special case of leaderless, the communication graph is not required to be connected.   When the initial states of followers  are identically and independently distributed, a near-optimal strategy is proposed for a large number of homogeneous followers  by solving two decoupled Riccati equations  using} mean-field team approach introduced by~\cite{arabneydi2016new} and showcased in~\cite{JalalCCECE2018, JalalCDC2017, JalalACC2018,JalalCDC2015}.

The remainder of this paper is organized as follows. The problem is formulated in the next section, where the main contributions of the work are also outlined. Then in Section~\ref{sec:main}, some important assumptions are presented and the control strategy is derived. Two numerical examples are provided in Section~\ref{sec:example} and finally the paper is concluded in Section~\ref{sec:conclusion}.

\section{Problem Formulation} \label{sec:prblmformul}
\subsection{Notation}
 In this paper, $\mathbb{N}$ and $\mathbb{R}$ denote natural and real numbers, respectively.  The short-hand notation $x_{a:b}$ is used to denote vector $(x_a, \ldots,x_b)$, $a \leq b \in \mathbb{N}$. For any $k \in \mathbb{N}$,  $\mathbb{N}_k$ is the finite set of integers  $\{1,2, \ldots,k\}$. $\TR(\boldsymbol \cdot)$ is the trace of a matrix and $\VAR(\boldsymbol \cdot)$ is the covariance of a random vector.

\subsection{Dynamics}
Consider a multi-agent system consisting  of   one leader and $n \in \mathbb{N}$ followers.  Denote by $x^0_t \in \mathbb{R}^{d_x},$ $u^0_t \in \mathbb{R}^{d_u}$, and $w^0_t \in \mathbb{R}^{d_x}$, $d_x,d_u \in \mathbb{N}$, the state, action, and noise of the leader at time $t \in \mathbb{N}$, respectively.  In addition, let $x^i_t \in \mathbb{R}^{d_x}$, $u^i_t \in \mathbb{R}^{d_u}$, and $w^i_t \in \mathbb{R}^{d_x}$ denote  the state, action, and noise of follower  $i \in \mathbb{N}_n$ at time $t \in \mathbb{N}$,  analogously.  The state of the leader  evolves as follows:
\begin{equation} \label{eq:dynamics-leader}
x^0_{t+1}=A_t^0x^0_t+B_t^0u^0_t+D^0_t \bar{x}_t+w^0_t,
\end{equation}
 where $\bar x_t:=\frac{1}{n} \sum_{i=1}^n x^i_t$ is the average of the states of the followers at time $t$, and will hereafter be called mean-field~\cite{arabneydi2016new}.  Similarly, the state of each follower $i \in \mathbb{N}_n$   evolves as follows:
\begin{equation} \label{eq:dynamics-followers}
x^i_{t+1}=A_tx^i_t+B_tu^i_t+D_t \bar{x}_t+E_t x_t^0+w^i_t.
\end{equation}
Let $\mathbb N_T $, $T \in \mathbb{N},$  be  the control horizon. It is assumed that the primitive random variables $\{x^0_1, x^1_1, \ldots, x^n_1, w^0_1,w^1_1,\ldots,w^n_1, \ldots, w^0_T,w^1_T,\ldots,w^n_T\}$ are defined on a common probability  space, and are mutually  independent.

\subsection{Information structure} \label{sec:info}
At time  $t \in \mathbb{N}$, the leader  observes its state $x^0_t$ and chooses its action $u^0_t$ according to a control law $g^0_t:\mathbb{R}^{d_x}  \rightarrow \mathbb{R}^{d_u}$, i.e.,
\begin{equation}   \label{eq:leader-inf. structre}
u^0_t=g^0_t(x^0_t).
\end{equation}
In addition, for any $i \in \mathbb{N}_n$, follower $i$  observes  its state $x^i_t$ as well as the state of the leader $x^0_t$ at time $t$, and decides  its action $u^i_t$  as follows:
\begin{equation} \label{eq:follower-inf. structre}
u^i_t=g^i_t(x^i_t,x^0_t),
\end{equation}
where $g^i_t:\mathbb{R}^{d_x} \times \mathbb{R}^{d_x} \rightarrow \mathbb{R}^{d_u}$ is the control law.
Note that the control actions~\eqref{eq:leader-inf. structre} and~\eqref{eq:follower-inf. structre} have a non-classical decentralized information structure.
\begin{remark}
The number of  links  required to implement  the above information structure  is $n$, which is the minimum possible number of links  for leader-follower networks to be connected. In addition, no information about the states of  the followers is communicated; hence,  the privacy of the  followers is preserved. 
\end{remark}

\begin{remark}
\edit{It is worth highlighting the difference between  the dynamics coupling~\eqref{eq:dynamics-leader} and~\eqref{eq:dynamics-followers}, that refers to the physical interactions among agents,  and  information coupling~\eqref{eq:leader-inf. structre} and~\eqref{eq:follower-inf. structre}, that attributes to the communication of data.  }
\end{remark}
The set  of  all control laws  $\mathbf g:=\{g^0_{1:T}, g^1_{1:T}, \ldots,g^n_{1:T}\}$ is called the strategy of the network. The objective of the followers is to  track the leader in an energy-efficient manner. To this end,  the following cost function is defined:
\begin{align}\label{eq:per-step}
&J_T(\mathbf g)\hspace{-.1cm}=\hspace{-.1cm}\mathbb{E} \Big[ \hspace{-.1cm} \sum_{t=1}^T  (x^0_t)^\intercal Q^0_t x^0_t \hspace{-.1cm}+\hspace{-.1cm}( \bar{x}_t \hspace{-.05cm}- \hspace{-.05cm}x^0_t)^\intercal F_t (\bar{x}_t \hspace{-.05cm}- \hspace{-.05cm}x^0_t) \hspace{-.1cm}+\hspace{-.1cm} (u^0_t)^\intercal R^0_t u^0_t \nonumber \\
& \quad +\frac{1}{n} \sum_{i=1}^n  (x^i_t)^\intercal Q_t x^i_t+( x^i_t \hspace{-.05cm}- \hspace{-.05cm}x^0_t)^\intercal P_t (x^i_t \hspace{-.05cm} - \hspace{-.05cm} x^0_t)+(u^i_t)^\intercal R_t u^i_t
\nonumber \\
&\quad +\frac{1}{2n^2} \sum_{i=1}^n   \sum_{j=1}^n  ( x^i_t-x^j_t)^\intercal H_t (x^i_t - x^j_t) \Big],
\end{align}
where the expectation is taken with respect to the probability measures induced by the choice of  strategy $\mathbf g$, and $\{Q^0_t, F_t, R^0_t, Q_t, P_t,  R_t, H_t\}$ are symmetric matrices of appropriate dimensions.  It is to be noted that the rate of convergence of the followers to the leader depends on the matrices $P_t$ and $F_t$.  Moreover, the collective behavior of the followers changes  by  matrix $H_t$.

\begin{Problem} \label{problem1}
	Consider the above leader-follower system with dynamics~\eqref{eq:dynamics-leader} and~\eqref{eq:dynamics-followers} and information structure~\eqref{eq:leader-inf. structre} and~\eqref{eq:follower-inf. structre}. We are interested to find an $\varepsilon(n)$-optimal strategy $\mathbf{g}^\ast_\varepsilon$ such that  for every strategy $\mathbf g$,
	\begin{equation}
	J_T(\mathbf{g}^\ast_\varepsilon) \leq J_T(\mathbf{g}) + \varepsilon(n),
	\end{equation}
	where  $\varepsilon(n) \in [0,\infty)$ and $\lim_{n \rightarrow \infty} \varepsilon(n)=0$.
\end{Problem}
\begin{remark}\label{remark:leaderless}
Notice that if matrices $B^0_t$, $D^0_t, Q^0_t$ and $R^0_t$ are zero,  Problem~\ref{problem1} reduces to the optimal control of a leaderless multi-agent network. In that case, $x^0_t$  represents the desired reference trajectory, and as noted before, the followers do not share anything with each other once they receive the reference trajectory information $x^0_t$, according to~\eqref{eq:follower-inf. structre}. 
\end{remark}
\subsection{Main challenges and contributions}
There are two main challenges in finding a solution to Problem~\ref{problem1}. The first one is concerned with non-classical information structure, as the optimal strategy under this type of information structure is not necessarily linear~\cite{Witsenhausen1968Counterexample}.  The second challenge is the curse of dimensionality as the matrices in Problem~\ref{problem1} are fully dense, yet  their dimension  increases with the number of followers.

In their previous work~\cite{JalalCCECE2018}, the authors show that if the mean-field $\bar{x}_t$ is available to the leader and followers, then the optimal solution is unique and linear.  However, collecting and sharing  the mean-field among all agents is not cost-efficient, in general, specially when  the number of followers $n$ is large. It is shown in the next section that the effect of such information sharing on the performance of the network is negligible when  the number of followers is large enough.

\section{Main Results}\label{sec:main}
In this section, we propose a strategy and compute its performance with respect to the optimal performance, and show that the difference between them converges to zero at rate $\frac{1}{n}$.  For the sake of clarity in the notation, we use letters $s$ and $v$ to denote the states and actions, respectively, under the optimal strategy.  Therefore, from~\eqref{eq:dynamics-leader} and ~\eqref{eq:dynamics-followers}, the dynamics of the leader and followers  at time $t \in \mathbb{N}_T$ under the optimal strategy are given by
\begin{equation} \label{eq:optimal dynamics-leader}
s^0_{t+1}=A_t^0s^0_t+B_t^0v^0_t+D^0_t \bar{s}_t+w^0_t,
\end{equation}
\begin{equation} \label{eq:optimal dynamics-followers}
s^i_{t+1}=A_ts^i_t+B_tv^i_t+D_t \bar{s}_t+E_t s_t^0+w^i_t, \quad i \in \mathbb{N}_n,
\end{equation}
where $\bar{s}_t:=\frac{1}{n} \sum_{i=1}^n s^i_t$. Similarly to~\cite{JalalCCECE2018}, define the following matrices:
\begin{equation}
\bar A_t \hspace{-.1cm}:= \hspace{-.1cm}\left[ \begin{array}{cc}
A^0_t &D^0_t\\
E_t &A_t+D_t
\end{array}  \right], \quad \bar B_t \hspace{-.1cm}:= \hspace{-.1cm}\left[ \begin{array}{cc}
B^0_t &\mathbf{0}_{d_x \times d_u}\\
\mathbf{0}_{d_x \times d_u} &B_t
\end{array}  \right],
\end{equation}
\begin{equation}
\bar Q_t \hspace{-.1cm}:= \hspace{-.1cm}\left[  \hspace{-.2cm}  \begin{array}{cc}
Q^0_t \hspace{-.1cm}+ \hspace{-.1cm}P_t \hspace{-.1cm}+ \hspace{-.1cm}F_t & - P_t \hspace{-.1cm}- \hspace{-.1cm}F_t\\
- P_t \hspace{-.1cm}- \hspace{-.1cm}F_t & Q_t  \hspace{-.1cm} +  \hspace{-.1cm}P_t+  \hspace{-.1cm}F_t
\end{array}  \hspace{-.2cm}   \right],  \bar R_t \hspace{-.1cm}:= \hspace{-.1cm}\left[\hspace{-.2cm}  \begin{array}{cc}
R^0_t &\mathbf{0}_{d_u \times d_u}\\
\mathbf{0}_{d_u \times d_u} &R_t
\end{array} \hspace{-.2cm} \right].
\end{equation} 
\begin{Assumption}\label{assump:positive}
Matrices  $Q_t + P_t + H_t$ and $\bar Q_t$  are positive semi-definite and  matrices $R_t$ and $R^0_t$ are positive definite.
\end{Assumption}
Define two decoupled  Riccati equations s.t. for  any $t \in \mathbb{N}_T$,
\begin{multline}\label{eq:Riccati-1}
\breve M_t = - A_t^\intercal   \breve M_{t+1} B_t  \left( B_t^\intercal  \breve M_{t+1}  B_t  +  R_t  \right)^{-1}  B_t^\intercal   \breve M_{t+1}  A_t\\
+ A_t^\intercal  \breve M_{t+1} A_t + Q_t+P_t+H_t,
\end{multline}
\begin{multline}\label{eq:Riccati-2}
\bar M_t = - \bar A_t^\intercal   \bar M_{t+1} \bar B_t  \left( \bar B_t^\intercal  \bar M_{t+1} \bar  B_t  +  \bar R_t  \right)^{-1}  \bar B_t^\intercal   \bar M_{t+1}  \bar A_t\\
+\bar  A_t^\intercal  \bar  M_{t+1} \bar A_t +\bar Q_t,
\end{multline}
where $\breve M_{T+1}=\mathbf{0}_{d_x \times d_x}$ and  $\bar M_{T+1}= \mathbf{0}_{2d_x \times 2d_x}$.  According to~\cite{arabneydi2016new,JalalCCECE2018}, the optimal performance $J_T^\ast$ is obtained under the following linear strategies:
\begin{align} \label{eq:optimal_action_leader}
{v^0_t}&=\bar L^{1,1}_t s^0_t +  \bar L^{1,2}_t \bar s_t, 
\end{align}
\begin{align} \label{eq:optimal_action_followers}
{v^i_t}&=\breve L_ts^i_t + \bar L_t^{2,1} s^0_t+ (\bar L_t^{2,2}-\breve L_t)\bar s_t,
\end{align}
where $\breve{L}_t$ and $\bar{L}_t=:\left[ 
\begin{array}{cc}
\bar L^{1,1}_t & \bar L^{1,2}_t\\
\bar L^{2,1}_t & \bar L^{2,2}_t
\end{array}
\right]$ can be found by using these formulas:
\begin{align}\label{eq:gains}
\breve L_t&=-\left( B_t^\intercal  \breve M_{t+1}  B_t  +  R_t  \right)^{-1}  B_t^\intercal   \breve M_{t+1}  A_t,\nonumber  \\
\bar L_t&=-\left( \bar B_t^\intercal  \bar M_{t+1}  \bar B_t  +  \bar R_t  \right)^{-1}  \bar B_t^\intercal   \bar M_{t+1}  \bar A_t.
\end{align}

\subsection{Solution of Problem~\ref{problem1}}
The following standard assumptions are imposed.
\begin{Assumption} \label{assump:noise}
	The  initial states of the followers are i.i.d.  with mean $\mu_x \in \mathbb{R}^{d_x}$  and  finite covariance matrix $\sum_x  \in \mathbb{R}^{d_x} \times  \mathbb{R}^{d_x}$.  In addition, the local noises of the followers are i.i.d. with zero mean and  finite covariance  matrix $\sum_w  \in \mathbb{R}^{d_x} \times  \mathbb{R}^{d_x}$.
\end{Assumption}

\begin{Assumption} \label{assmpt:variance bounded}
 Matrices    $A_t, A_t^0, B_t, B_t^0, D_t, D_t^0, E_t, F_t, Q_t$, $ Q_t^0$, $ R_t$,  $R_t^0$, $\sum_x$ and $\sum_w$   do not depend on the number of followers $n$.
\end{Assumption} 
Define a stochastic process $z_{1:T}$ such that $z_1:=\mu_x$ and for any $t \in \mathbb{N}_T$:
\begin{equation} \label{eq:approximation process}
z_{t+1}=(A_t+B_t\bar{L}_t^{2,2}+D_t)z_t+(B_t\bar{L}_t^{2,1}+E_t)x_t^0.
\end{equation}
Note that the  leader and followers can compute $z_t$ under the information structures~\eqref{eq:leader-inf. structre} and~\eqref{eq:follower-inf. structre}.  Given the matrix gains defined in~\eqref{eq:gains}, the following strategies are proposed:
\begin{equation} \label{eq:our strategy for leader}
u_t^0=\bar{L}_t^{1,1} x_t^0 + \bar{L}^{1,2}_t z_t,
\end{equation}
\begin{equation} \label{eq:our strategy for followers}
u_t^i=\breve{L}_t x_t^i + \bar{L}^{2,1}_t x_t^0 + (\bar{L}^{2,2}_t - \breve L_t) z_t, \quad i \in \mathbb{N}_n.
\end{equation}
At any time $t \in \mathbb{N}_T$, define the following relative errors  $e_t^0$, $e_t$ and  $\zeta_t$:
\begin{equation}\label{eq:relative-error}
e_t^0:=s_t^0-x_t^0, \quad e_t:=\bar{s}_t-z_t, \quad \zeta_t:=\bar{x}_t-z_t.
\end{equation}  
    \begin{Lemma} \label{lemma:error dynamics}
	The relative errors defined in~\eqref{eq:relative-error} evolve linearly in time as follows:  
	\begin{equation} \label{eq:final error equation}
	\left[ \begin{array}{c}
	e_{t+1}^0\\
	e_{t+1} \\
	\zeta_{t+1}  
	\end{array}  \right]=
	\tilde{A}_t  \left[ \begin{array}{c}
	e_{t}^0\\
	e_{t} \\
	\zeta_{t}
	\end{array}  \right]+\left[ \begin{array}{c}
	\mathbf{0}_{d_x \times 1}\\
	\bar{w}_{t}\\  
	\bar{w}_{t} 
	\end{array}  \right],
	\end{equation}
	where $\bar w_t:=\frac{1}{n}\sum_{i=1}^n w^i_t$ and
	\begin{equation}
	\tilde{A}_t\hspace{-.1cm} := \hspace{-.1cm} \left[ \hspace{-.2cm} \begin{array}{c c c}
	A_t^0+B_t^0 \bar{L}_t^{1,1} & B_t^0 \bar{L}_t^{1,2}+ D_t^0 & -D_t^0 \\
	B_t \bar{L}_t^{2,1} +E_t & A_t+B_t \bar{L}_t^{2,2}+D_t & \mathbf{0}_{d_x \times d_x} \\
	\mathbf{0}_{d_x \times d_x} & \mathbf{0}_{d_x \times d_x} & A_t+ B_t \breve{L}_t +D_t
	\end{array} \hspace{-.2cm}  \right].
	\end{equation}
\end{Lemma}

\begin{proof}
	From~\eqref{eq:optimal dynamics-leader} and~\eqref{eq:optimal_action_leader}:
	\begin{equation} \label{eq:s_0_bar_dynamics after applying control}
	s_{t+1}^0=(A_t^0 + B_t^0 \bar{L}^{1,1}_t)s_t^0 + (B_t^0 \bar{L}_t^{1,2}+D_t^0) \bar{s}_t+w_t^0. 
	\end{equation}
	Also, it results from \eqref{eq:dynamics-leader} and~\eqref{eq:our strategy for leader} that
	\begin{equation} \label{eq:x_0_bar_dynamics after applying control}
	x_{t+1}^0=(A_t^0 + B_t^0 \bar{L}^{1,1}_t)x_t^0 + B_t^0 \bar{L}_t^{1,2} z_t + D_t^0 \bar{x}_t+w_t^0. 
	\end{equation}
	Similarly, from~\eqref{eq:optimal dynamics-followers} and~\eqref{eq:optimal_action_followers}:
	\begin{equation} \label{eq: average_optimal_dynamics} 
	\bar{s}_{t+1}=A_t\bar{s}_t+B_t\bar{v}_t+D_t \bar{s}_t+E_t s_t^0+\bar{w}_t,
	\end{equation}
	where $\bar{v}_t:=\frac{1}{n}\sum_{i=1}^n v^i_t$ is given by: 
	\begin{equation} \label{eq:average_our_strategy}
	\bar{v}_t=\breve L_t \bar{s}_t + \bar L_t^{2,1} s^0_t+ (\bar L_t^{2,2}-\breve L_t)\bar s_t
	=\bar L_t^{2,1} s^0_t+\bar L_t^{2,2} \bar s_t.
	\end{equation}
	Substituting ~\eqref{eq:average_our_strategy} in ~\eqref{eq: average_optimal_dynamics} yields:
	\begin{equation}\label{eq:s_bar_dynamics with control}
	\bar{s}_{t+1}=(A_t+B_t \bar{L}_t^{2,2}+D_t) \bar{s}_t+(B_t \bar{L}_t^{2,1}+E_t) s_t^0+\bar{w}_t.
	\end{equation}
	In addition, from~\eqref{eq:dynamics-followers} and~\eqref{eq:our strategy for followers}, one arrives at:
		\begin{equation} \label{eq:x_bar_dynamics without control}
		\bar{x}_{t+1}=A_t \bar{x}_t+B_t \bar{u}_t+D_t \bar{x}_t+E_t x_t^0+\bar{w}_t,
		\end{equation}
		where $\bar{u}_t:=\frac{1}{n}\sum_{i=1}^n u^i_t$ is as follows:
		\begin{equation} \label{eq:x_bar_control}
		\bar{u}_{t+1}=\breve{L}_t \bar{x}_t+ \bar{L}_t^{2,1}x_t^0+ (\bar{L}_t^{2,2} - \breve L_t)z_t.
		\end{equation}
		From~\eqref{eq:x_bar_dynamics without control} and~\eqref{eq:x_bar_control}, it results that:
		\begin{multline} \label{eq:x_bar_dynamics with control}
		\bar{x}_{t+1}=A_t \bar{x}_t+B_t \breve{L}_t \bar{x}_t+B_t \bar{L}_t^{2,1}x_t^0+B_t (\bar{L}_t^{2,2} - \breve L_t)z_t
		\\+D_t \bar{x}_t+E_t x_t^0+\bar{w}_t.
		\end{multline}
	Equations~\eqref{eq:relative-error},~\eqref{eq:s_0_bar_dynamics after applying control} and~\eqref{eq:x_0_bar_dynamics after applying control} lead to:	
	\begin{align}	\label{eq:e_t^0}
	e_{t+1}^0&=(A_t^0 +B_t^0 \bar{L}_t^{1,1})s_t^0+B_t^0 \bar{L}_t^{1,2} \bar{s}_t+ D_t^0 \bar{s}_t+w_t^0
	\nonumber 	\\
	&\quad -(A_t^0 x_t^0 +B_t^0 \bar{L}_t^{1,1}x_t^0+B_t^0 \bar{L}_t^{1,2} z_t+D_t^0 \bar{x}_t+w_t^0)\nonumber \\
	&=(A_t^0 + B_t^0 \bar{L}_t^{1,1}) e_t^0 + (B_t^0 \bar{L}_t^{1,2}+ D_t^0) e_t -D_t^0 \zeta_{t}.
	\end{align}
	Moreover, it results from ~\eqref{eq:approximation process},~\eqref{eq:relative-error} and~\eqref{eq:s_bar_dynamics with control} that:
	\begin{align}	\label{eq:e}
	e_{t+1}&= A_t \bar{s}_t + B_t \bar{L}^{2,2}_t \bar{s}_t +D_t \bar{s}_t + B_t \bar{L}_t^{2,1}s_t^0 + E_t s_t^0 + \bar{w}_t
	\nonumber \\&\quad -(A_t z_t+B_t\bar{L}_t^{2,2} z_t+D_t z_t + B_t\bar{L}_t^{2,1} x_t^0+E_t x_t^0),\\
	&=(A_t +B_t \bar{L}_t^{2,2} + D_t) e_{t}+(B_t \bar{L}_t^{2,1} +E_t) e_{t}^0 +\bar{w}_t.
	\end{align}
	As a consequence of~\eqref{eq:approximation process},~\eqref{eq:relative-error} and~\eqref{eq:x_bar_dynamics with control}, the following equation is obtained: 
	\begin{align}	 \label{eq:zeta}
	&\zeta_{t+1}= (A_t+B_t \breve{L}_t+D_t) \bar{x}_t  +B_t (\bar{L}_t^{2,1} - \breve L_t) x_t^0 +B_t \bar{L}_t^{2,2} z_t \nonumber \\
	& +E_t x_t^0 +\bar{w}_t
	-(A_t +B_t\bar{L}_t^{2,2}+D_t )z_t -( B_t\bar{L}_t^{2,1}+E_t)x_t^0 \nonumber \\
	& \qquad =(A_t + B_t \breve{L}_t+D_t) \zeta_{t} +\bar{w}_t,
	\end{align}
	and this completes the proof. $\hfill \blacksquare$
\end{proof}  

Now, for any follower $i \in \mathbb{N}_n$, define the following variables at time  $t \in \mathbb{N}_T$:
\begin{equation}\label{eq:def-breve}
\breve{x}_t^i:=x_t^i-\bar{x}_t,  \breve{u}_t^i:=u_t^i-\bar{u}_t, \breve{s}_t^i:=s_t^i-\bar{s}_t, \breve{v}_t^i:=v_t^i-\bar{v}_t.
\end{equation}
\begin{Lemma} \label{breve x and breve s}
	At any time $t \in \mathbb{N}_T$, $\breve{x}_t^i =\breve{s}_t^i $ and $\breve{u}_t^i=\breve{v}_t^i$.
\end{Lemma}
\begin{proof}
	The lemma is proved by induction on noting that initially $\breve x^i_1=\breve s^i_1=x^i_1- \bar x_1$ because $x^i_1=s^i_1$. It follows from~\eqref{eq:optimal_action_followers} and~\eqref{eq:our strategy for followers} that $\breve{u}_1^i=\breve{v}_1^i=\breve L_1 (x^i_1- \bar x_1)$. 
	Suppose  $\breve{x}_t^i =\breve{s}_t^i$ and $\breve{u}_t^i =\breve{v}_t^i$ at time $t$. It is now desired to  show that $\breve{x}_{t+1}^i =\breve{s}_{t+1}^i$ and $\breve{u}_{t+1}^i=\breve{v}_{t+1}^i $. From~\eqref{eq:dynamics-followers} and~\eqref{eq:optimal dynamics-followers} and the induction assumption at $t=1$, one arrives at:
	\begin{equation}\label{eq:proof-breve}
	\breve{s}_{t+1}^i = A_t \breve{s}_{t}^i + B_t \breve{v}_t^i + \breve{w}_t^i 
	=A_t \breve{x}_{t}^i + B_t \breve{u}_t^i + \breve{w}_t^i = \breve{x}_{t+1}^i,
	\end{equation}
	where $\breve w^i_t:=w^i_t -\bar w_t$. Also, it is implied from~\eqref{eq:optimal_action_followers},~\eqref{eq:our strategy for followers} and~\eqref{eq:proof-breve} that:
	$\breve{v}_{t+1}^i =\breve{L}_{t+1} \breve{s}_{t+1}^i =\breve{L}_{t+1} \breve{x}_{t+1}^i = \breve{u}_{t+1}^i$. 
	$\hfill \blacksquare$
\end{proof}
\begin{Lemma} \label{Lemma: Error cost function}
	Let $\Delta J$ denote the discrepancy between the performance of the optimal strategies~\eqref{eq:optimal_action_leader} and~\eqref{eq:optimal_action_followers}, and that of the proposed  strategies~\eqref{eq:our strategy for leader} and~\eqref{eq:our strategy for followers}. If Assumption~\ref{assump:noise} holds,  $\Delta J$ is a quadratic function of the relative errors in~\eqref{eq:relative-error}, i.e.,
	\begin{equation}
	\Delta J=\mathbb{E} \Big[\sum_{t=1}^T  \left[\begin{array}{c c c}
	e_t ^0 & e_t & \zeta_{t} 
	\end{array}\right]^\intercal \tilde{Q}_t \left[\begin{array}{c c c}
	e_t ^0 & e_t & \zeta_{t} 
	\end{array}\right] \Big],
	\end{equation}
	where
	\begin{equation}
	\tilde{Q}_t:=\left[\begin{array}{c c}
	-\bar{Q}_t - \bar{L}_t^\intercal \bar{R}_t \bar{L}_t & \mathbf{0}_{2d_x \times d_x}\\
	\mathbf{0}_{d_x \times 2d_x} & Q_t+P_t+F_t+\breve{L}_t ^\intercal R_t \breve{L}_t
	\end{array}\right].
	\end{equation}
\end{Lemma}
\begin{proof}
	From~\eqref{eq:per-step}, we have
	\begin{multline}
	\Delta J=\mathbb{E} \Big[\sum_{t=1}^T  (x^0_t)^\intercal Q^0_t x^0_t + (u^0_t)^\intercal R^0_t u^0_t+(\bar{x}_t - x^0_t)^\intercal F_t (\bar{x}_t - x^0_t)\\
	+\frac{1}{n} \sum_{i=1}^n  (x^i_t)^\intercal Q_t x^i_t+( x^i_t-x^0_t)^\intercal P_t (x^i_t - x^0_t)
	+ (u^i_t)^\intercal R_t u^i_t \\
	+\frac{1}{2n^2} \sum_{i=1}^n   \sum_{j=1}^n  ( x^i_t-x^j_t)^\intercal H_t (x^i_t - x^j_t) \Big]
	\\-\mathbb{E} \Big[\sum_{t=1}^T  (s^0_t)^\intercal Q^0_t s^0_t + (v^0_t)^\intercal R^0_t v^0_t+(\bar{s}_t - s^0_t)^\intercal F_t (\bar{s}_t - s^0_t)\\
	+\frac{1}{n} \sum_{i=1}^n  (s^i_t)^\intercal Q_t s^i_t+( s^i_t-s^0_t)^\intercal P_t (s^i_t - s^0_t)+ (v^i_t)^\intercal R_t v^i_t\\
	+\frac{1}{2n^2} \sum_{i=1}^n   \sum_{j=1}^n  ( s^i_t-s^j_t)^\intercal H_t (s^i_t - s^j_t) \Big].
	\end{multline}
	The above equation can be re-written in terms of the variables defined in~\eqref{eq:def-breve} as follows:
	\begin{multline}\label{eq:proof_sentence1}
	\Delta J=\mathbb{E} \Big[\sum_{t=1}^T  (x^0_t)^\intercal Q^0_t x^0_t + (u^0_t)^\intercal R^0_t u^0_t
	+(\bar{x}_t - x^0_t)^\intercal F_t (\bar{x}_t - x^0_t)
	\\
	+\frac{1}{n} \sum_{i=1}^n  (\breve{x}^i_t+\bar{x}_t)^\intercal Q_t (\breve{x}^i_t+\bar{x}_t) +(\breve{x}^i_t+\bar{x}_t - x^0_t)^\intercal P_t (\breve{x}^i_t+\bar{x}_t - x^0_t)
	\\
	+(\breve{u}_t^i+ \bar{u}_t)^\intercal R_t (\breve{u}_t^i+ \bar{u}_t)	
	+\frac{1}{2n^2} \sum_{i=1}^n   \sum_{j=1}^n  ( \breve{x}^i_t-\breve{x}^j_t)^\intercal H_t (\breve{x}^i_t - \breve{x}^j_t) \Big]
	\\
	-\mathbb{E} \Big[\sum_{t=1}^T  (s^0_t)^\intercal Q^0_t s^0_t + (v^0_t)^\intercal R^0_t v^0_t+(\bar{s}_t - s^0_t)^\intercal F_t (\bar{s}_t - s^0_t)\\
	+\frac{1}{n} \sum_{i=1}^n  (\breve{s}^i_t+ \bar{s}_t)^\intercal Q_t (\breve{s}^i_t+ \bar{s}_t)+(\breve{s}^i_t+ \bar{s}_t-s^0_t)^\intercal P_t (\breve{s}^i_t+ \bar{s}_t - s^0_t)\\+(\breve{v}^i_t +\bar{v}_t)^\intercal R_t (\breve{v}^i_t +\bar{v}_t)
	+\frac{1}{2n^2} \sum_{i=1}^n   \sum_{j=1}^n  ( \breve{s}^i_t-\breve{s}^j_t)^\intercal H_t (\breve{s}^i_t - \breve{s}^j_t) \Big].
	\end{multline} 
	On the other hand, by definition the following relations hold:
	\begin{equation}\label{eq:relations}
	\sum_{i=1}^{n} \breve{x}_t^i=\sum_{i=1}^{n} \breve{s}_t^i=\mathbf{0}_{d_x \times 1},\sum_{i=1}^{n} \breve{u}_t^i=\sum_{i=1}^{n} \breve{v}_t^i=\mathbf{0}_{d_u \times 1}.
	\end{equation}
	By incorporating~\eqref{eq:relations} in~\eqref{eq:proof_sentence1}, it results that
	\begin{multline}
	\Delta J\hspace{-.1cm} = \hspace{-.1cm}\mathbb{E} \Big[\sum_{t=1}^T  \left[\begin{array}{c}
	x^0_t\\
	\bar x_t
	\end{array}\right]^\intercal \bar Q_t \left[\begin{array}{c}
	x^0_t\\
	\bar x_t
	\end{array}\right] \hspace{-.1cm}+ \hspace{-.1cm}  \left[\begin{array}{c}
	u^0_t\\
	\bar u_t
	\end{array}\right]^\intercal \bar R_t  \left[\begin{array}{c}
	u^0_t\\
	\bar u_t
	\end{array}\right] \nonumber \\
	+\frac{1}{n}\sum_{i=1}^n (\breve  x^i_t)^\intercal (Q_t+P_t+H_t) (\breve x^i_t)+ (\breve u^i_t)^\intercal R_t (\breve u^i_t) \Big]
	\\ -\mathbb{E} \Big[\sum_{t=1}^T  \left[\begin{array}{c}
	s^0_t\\
	\bar s_t
	\end{array}\right]^\intercal \bar Q_t \left[\begin{array}{c}
	s^0_t\\
	\bar s_t
	\end{array}\right] \hspace{-.1cm} + \hspace{-.1cm}  \left[\begin{array}{c}
	v^0_t\\
	\bar v_t
	\end{array}\right]^\intercal \bar R_t  \left[\begin{array}{c}
	v^0_t\\
	\bar v_t
	\end{array}\right] \nonumber \\
	+\frac{1}{n}\sum_{i=1}^n (\breve  s^i_t)^\intercal (Q_t+P_t+H_t) (\breve s^i_t)+ (\breve v^i_t)^\intercal R_t (\breve v^i_t) \Big].
	\end{multline}
	According to Lemma~\ref{breve x and breve s}, the above equation simplifies to
	\begin{multline} \label{eq:decoupled form}
	\Delta J \hspace{-.1cm}= \hspace{-.1cm} \mathbb{E} \Big[\sum_{t=1}^T  \left[\begin{array}{c}
	x^0_t\\
	\bar x_t
	\end{array}\right]^\intercal \bar Q_t \left[\begin{array}{c}
	x^0_t\\
	\bar x_t
	\end{array}\right]+  \left[\begin{array}{c}
	u^0_t\\
	\bar u_t
	\end{array}\right]^\intercal \bar R_t  \left[\begin{array}{c}
	u^0_t\\
	\bar u_t
	\end{array}\right]\Big]
	\\ -\mathbb{E} \Big[\sum_{t=1}^T  \left[\begin{array}{c}
	s^0_t\\
	\bar s_t
	\end{array}\right]^\intercal \bar Q_t \left[\begin{array}{c}
	s^0_t\\
	\bar s_t
	\end{array}\right]+  \left[\begin{array}{c}
	v^0_t\\
	\bar v_t
	\end{array}\right]^\intercal \bar R_t  \left[\begin{array}{c}
	v^0_t\\
	\bar v_t
	\end{array}\right] \Big].
	\end{multline}
	Using the definition of relative errors in~\eqref{eq:relative-error}, one concludes that:
	\begin{multline}\label{eq:proof_sentence2}
	\hspace{-.4cm}\left[\begin{array}{c}
	x^0_t\\
	\bar x_t
	\end{array}\right] = \left[\begin{array}{c}
	x^0_t\\
	z_t
	\end{array}\right] + \left[\begin{array}{c}
	\mathbf{0}_{d_x \times d_x}\\
	\zeta_{t}
	\end{array}\right], 	
	\left[\begin{array}{c}
	s^0_t\\
	\bar s_t
	\end{array}\right] =\left[\begin{array}{c}
	e^0_t+x_t^0\\
	e_t+z_t
	\end{array}\right], \\
	\left[\begin{array}{c}
	u^0_t\\
	\bar u_t
	\end{array}\right] = \bar{L}_t \left[\begin{array}{c}
	x_t^0\\
	z_t
	\end{array}\right]+\left[\begin{array}{c}
	\mathbf{0}_{d_x \times d_x}\\
	\breve{L}_t  \zeta_{t}
	\end{array}\right],  \left[\begin{array}{c}
	v^0_t\\
	\bar v_t
	\end{array}\right] =  \bar{L}_t \left[\begin{array}{c}
	s^0_t\\
	\bar{s}_t
	\end{array}\right].
	\end{multline}
	It is implied from~\eqref{eq:decoupled form} and~\eqref{eq:proof_sentence2} that:
	\begin{multline}\label{eq:proof_sentence3} 
	\Delta J\hspace{-.1cm} = \hspace{-.1cm} \mathbb{E} \Big[\sum_{t=1}^T (\left[ \hspace{-.1cm} \begin{array}{c}
	x_t^0\\
	z_t
	\end{array} \hspace{-.1cm} \right] \hspace{-.05cm} + \hspace{-.05cm} \left[\hspace{-.1cm} \begin{array}{c}
	\mathbf{0}_{d_x \times d_x}\\
	\zeta_{t}
	\end{array} \hspace{-.1cm} \right])^\intercal \bar Q_t (\left[ \hspace{-.1cm} \begin{array}{c}
	x_t^0\\
	z_t
	\end{array} \hspace{-.1cm} \right]+\left[\hspace{-.1cm} \begin{array}{c}
	\mathbf{0}_{d_x \times d_x}\\
	\zeta_{t}
	\end{array} \hspace{-.1cm} \right])\Big]
	\\
	+\mathbb{E} \Big[\sum_{t=1}^T( \bar{L}_t \left[ \hspace{-.1cm}  \begin{array}{c}
	x_t^0\\
	z_t
	\end{array} \hspace{-.1cm}  \right]+\left[ \hspace{-.1cm}  \begin{array}{c}
	\mathbf{0}_{d_x \times d_x}\\
	\breve{L}_t  \zeta_{t}
	\end{array} \hspace{-.1cm}  \right] )^\intercal \bar R_t  (\bar{L}_t \left[\hspace{-.1cm}  \begin{array}{c}
	x_t^0\\
	z_t
	\end{array} \hspace{-.1cm}  \right]\hspace{-.05cm} +\hspace{-.05cm} \left[\hspace{-.1cm}  \begin{array}{c}
	\mathbf{0}_{d_x \times d_x}\\
	\breve{L}_t \zeta_{t}
	\end{array} \hspace{-.1cm} \right] )\Big]
	\\
	-\mathbb{E} \Big[\sum_{t=1}^T  \left[\begin{array}{c}
	e^0_t+x_t^0\\
	e_t+z_t
	\end{array}\right]^\intercal (\bar Q_t+\bar{L}_t^\intercal \bar R_t \bar{L}_t ) \left[\begin{array}{c}
	e^0_t+x_t^0\\
	e_t+z_t
	\end{array}\right]\Big].
	\end{multline}
	Expand~\eqref{eq:proof_sentence3} as follows:
	\begin{align} \label{eq:proof_sentence4}
	\Delta J=&  \mathbb{E} \Big[\sum_{t=1}^T \left[\begin{array}{c}
	x_t^0\\
	z_t
	\end{array}\right]^\intercal \bar Q_t \left[\begin{array}{c}
	x_t^0\\
	z_t
	\end{array}\right]\Big]
	\nonumber \\ 
	&+2\mathbb{E} \Big[\sum_{t=1}^T \left[\begin{array}{c}
	x_t^0\\
	z_t
	\end{array}\right]^\intercal \bar Q_t  \left[\begin{array}{c}
	\mathbf{0}_{d_x \times d_x}\\
	\zeta_{t}
	\end{array}\right] \Big]
	\nonumber \\  
	&+\mathbb{E} \Big[\sum_{t=1}^T  \left[\begin{array}{c}
	\mathbf{0}_{d_x \times d_x}\\
	\zeta_{t}
	\end{array}\right]^\intercal  \bar Q_t \left[\begin{array}{c}
	\mathbf{0}_{d_x \times d_x}\\
	\zeta_{t}
	\end{array}\right] \Big]
	\nonumber \\ 
	&+\mathbb{E} \Big[\sum_{t=1}^T \left[\begin{array}{c}
	x_t^0\\
	z_t
	\end{array}\right]^\intercal \bar{L}_t ^\intercal \bar R_t  \bar{L}_t \left[\begin{array}{c}
	x_t^0\\
	z_t
	\end{array}\right] \Big]
	\nonumber \\ 
	&+2\mathbb{E} \Big[\sum_{t=1}^T   \left[\begin{array}{c}
	\mathbf{0}_{d_x \times d_x}\\
	\breve{L}_t  \zeta_{t}
	\end{array}\right]^\intercal \bar R_t \bar{L}_t \left[\begin{array}{c}
	x_t^0\\
	z_t
	\end{array}\right]\Big]
	\nonumber \\
	&+\mathbb{E} \Big[\sum_{t=1}^T \left[\begin{array}{c}
	\mathbf{0}_{d_x \times d_x}\\
	\breve{L}_t  \zeta_{t}
	\end{array}\right]^\intercal \bar R_t \left[\begin{array}{c}
	\mathbf{0}_{d_x \times d_x}\\
	\breve{L}_t  \zeta_{t}
	\end{array}\right]\Big]
	\nonumber \\ 
	&-\mathbb{E} \Big[\sum_{t=1}^T \left[\begin{array}{c}
	e_t^0\\
	e_t
	\end{array}\right]^\intercal (\bar Q_t +\bar{L}_t^\intercal \bar R_t \bar{L}_t )\left[\begin{array}{c}
	e_t^0\\
	e_t
	\end{array}\right]\Big]
	\nonumber \\ 
	&-2\mathbb{E} \Big[\sum_{t=1}^T \left[\begin{array}{c}
	e_t^0\\
	e_t
	\end{array}\right]^\intercal (\bar Q_t+\bar{L}_t^\intercal \bar R_t \bar{L}_t )  \left[\begin{array}{c}
	x_t^0\\
	z_t
	\end{array}\right] \Big]
	\nonumber \\ 
	&-\mathbb{E} \Big[\sum_{t=1}^T \left[\begin{array}{c}
	x_t^0\\
	z_t
	\end{array}\right]^\intercal (\bar Q_t+\bar{L}_t^\intercal \bar R_t \bar{L}_t )  \left[\begin{array}{c}
	x_t^0\\
	z_t
	\end{array}\right] \Big].
	\end{align}
	The second, fifth and eighth terms in the right side of~\eqref{eq:proof_sentence4} are zero from Lemma~\ref{lemma:error dynamics} and Assumption 2, on noting that $x^0_t$ and $z_t$ are completely  known under the information structures~\eqref{eq:leader-inf. structre} and~\eqref{eq:follower-inf. structre}. This completes the proof. $\hfill \blacksquare$
\end{proof}
\begin{Theorem}
	Let Assumptions~\ref{assump:positive} and \ref{assump:noise} hold. Then, 
	\begin{align} \label{eq: delta J}
	\Delta J&=\TR\left(    
	\left[ \begin{array}{ccc}
	\mathbf{0}_{d_x \times d_x} & \mathbf{0}_{d_x \times d_x} & \mathbf{0}_{d_x \times d_x}\\
	\mathbf{0}_{d_x \times d_x} & \VAR(\bar x_1)  & \VAR(\bar x_1)\\
	\mathbf{0}_{d_x \times d_x} &  \VAR(\bar x_1) & \VAR(\bar x_1)
	\end{array}
	\right] \tilde{M}_1 \right) \nonumber \\
	&\quad +\sum_{t=1}^{T-1}\TR\left(
	\left[ \begin{array}{ccc}
	\mathbf{0}_{d_x \times d_x} & \mathbf{0}_{d_x \times d_x} & \mathbf{0}_{d_x \times d_x}\\
	\mathbf{0}_{d_x \times d_x} & \VAR(\bar w_t)  & \VAR(\bar w_t)\\
	\mathbf{0}_{d_x \times d_x} &  \VAR(\bar w_t) & \VAR(\bar w_t)
	\end{array}
	\right] 
	\tilde{M}_{t+1}\right),
	\end{align}
	where  $\tilde{M}_{T}=\tilde{Q}_T$,  and $\tilde{M}_t$ is the solution of the following  Lyapunov equation for any  $t \in \mathbb{N}_{T-1}$ \textnormal{:}
	\begin{equation}\label{eq:lyapunov}
	\tilde{M}_t=\tilde{A}_t^T \tilde{M}_{t+1} \tilde{A}_t+\tilde{Q}_t.
	\end{equation} 
	
\end{Theorem}
\begin{proof}
	According to Lemma~\ref{Lemma: Error cost function}, the performance discrepancy $\Delta J$ is a quadratic function of the relative errors, and from Lemma~\ref{lemma:error dynamics}, the relative errors  have linear dynamics. Therefore, $\Delta J$ can be regarded as the quadratic cost of an uncontrolled linear system (where there is no control action). Thus, from the  standard results in linear systems~\cite{caines1987linear}, $\Delta J$ can be expressed by the Lyaponuv equation~\eqref{eq:lyapunov} and the covariance matrices of the initial relative errors and noises $\bar w_t$, $t \in \mathbb N_{T-1}$ as:
	\begin{multline}
	\Exp{[e^0_1 \quad e_1-\mu_x \quad \zeta_1-\mu_x]^\intercal [e^0_1 \quad e_1-\mu_x \quad \zeta_1-\mu_x]}\\=\left[ \begin{array}{ccc}
	\mathbf{0}_{d_x \times d_x} & \mathbf{0}_{d_x \times d_x} & \mathbf{0}_{d_x \times d_x}\\
	\mathbf{0}_{d_x \times d_x} & \VAR(\bar x_1) & \VAR(\bar x_1)\\
	\mathbf{0}_{d_x \times d_x} &  \VAR(\bar x_1) & \VAR(\bar x_1)
	\end{array}
	\right],
	\end{multline}
	and 
	\begin{multline}
	\Exp{[\mathbf{0}_{d_x \times d_x} \quad \bar w_t \quad \bar w_t]^\intercal [\mathbf{0}_{d_x \times d_x} \quad \bar w_t \quad \bar w_t]}\\=\left[ \begin{array}{ccc}
	\mathbf{0}_{d_x \times d_x} & \mathbf{0}_{d_x \times d_x} & \mathbf{0}_{d_x \times d_x}\\
	\mathbf{0}_{d_x \times d_x} & \VAR(\bar w_t) & \VAR(\bar w_t)\\
	\mathbf{0}_{d_x \times d_x} &  \VAR(\bar w_t) &  \VAR(\bar w_t)
	\end{array}
	\right].
	\end{multline}	
\end{proof}

\begin{Theorem}\label{theorem2}
	Let  Assumptions~\ref{assump:positive},~\ref{assump:noise} and~\ref{assmpt:variance bounded} hold. Then, the strategies proposed in~\eqref{eq:our strategy for leader} and~\eqref{eq:our strategy for followers} are $\varepsilon(n)$-optimal solutions for Problem~\ref{problem1} such that
	\begin{equation}
	|J_T(g_\epsilon^*)-J_T^*| \leq \varepsilon(n) \in \mathcal{O}(\frac{1}{n}).
	\end{equation}
\end{Theorem}

\begin{proof}
	According to Assumption~\ref{assump:noise}, 
	\begin{align}\label{eq:proof2}
	\VAR(\bar x_1)= \VAR(\frac{1}{n} \sum_{i=1}^n x^i_1)=\frac{n \sum_x}{n^2}  =\frac{ \sum_x}{n}, \nonumber \\
	\VAR(\bar w_t)= \VAR(\frac{1}{n} \sum_{i=1}^n w^i_t)=\frac{n \sum_w}{n^2}  =\frac{ \sum_w}{n}.
	\end{align}
	In addition, from Assumption~\ref{assmpt:variance bounded}, matrices $\tilde{A}_t$ and  $\tilde{Q}_t$ given by Lemmas~\ref{lemma:error dynamics} and~\ref{Lemma: Error cost function} are independent  of the number of followers $n$, and so is $\tilde{M}_t$. Therefore, the performance discrepancy in~\eqref{eq: delta J} converges to zero at rate~$\mathcal{O}(\frac{1}{n})$ according to~\eqref{eq:proof2}. $\hfill \blacksquare$
\end{proof}
\begin{corollary}
		For the special case of a leaderless multi-agent network, let $x_{t+1}^0=x_t^0=\bar{x}_1$, $t \in \mathbb{N}_{T}$.
		Then, according to Theorem~\ref{theorem2}, strategy~\eqref{eq:our strategy for followers} steers all the followers to the initial mean $\bar{x}_1$ as $n$ grows to infinity. In addition, if the initial mean $\bar{x}_1$ is not known, it can be replaced by its expectation, i.e., $x_{t+1}^0 = x_t^0= \mu_x$, $t \in \mathbb{N}_{T}$, and the resultant strategy~\eqref{eq:our strategy for followers} steers all the followers to the initial mean consensus as $n$ grows to infinity, due to the strong law of large numbers.
\end{corollary}
\section{Numerical Examples}\label{sec:example}
\begin{figure}[h!]
	\centering
	\vspace{-3.8cm}
	\includegraphics[width=\linewidth]{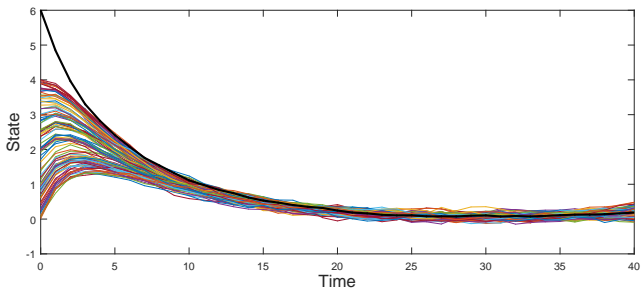} 
	\vspace{-4.2cm}
	\caption{Trajectories of the leader and $100$ randomly selected  followers in Example 1.}\label{fig 1}
\end{figure}

\textbf{Example 1:} Consider a multi-agent network with one leader and $1000$ followers, where the initial state of  the leader is $x^0_1=6$ and the initial states of the followers are chosen as uniformly distributed random variables in the interval $[0,4]$.  Let  the dynamics of the leader and followers  be described by~\eqref{eq:dynamics-leader} and \eqref{eq:dynamics-followers}, respectively, where 
\begin{align}
&\quad A^0_t=1, \quad B^0_t=0.8, \quad A_t=1, \quad  B_t=0.9, \\ 
&\quad D^0_t=0.1, \quad D_t=0.05, \quad E_t=0.15, \quad T=40, \\  
& \quad w^0_t \sim \mathcal{N}(0,0.02), \quad w^i_t \sim \mathcal{N}(0,0.05) \quad \forall i \in \mathbb{N}_{1000}.
\end{align}
The network objective is to minimize the cost function~\eqref{eq:per-step}, where
$ Q^0_t=1,  R^0_t=200,  F_t=20, Q_t=2,    P_t=5,   R_t=100,  H_t=1$. The leader  solves the Riccati equation~\eqref{eq:Riccati-2} to obtain gains $\bar L^{1,1}_t$ and $\bar L^{1,2}_t$, $t \in \mathbb{N}_T,$ and determines its control action according to strategy~\eqref{eq:our strategy for leader} using its local state $x^0_t$ as well as $z_{t}$. It is to be noted that  $z_{t}$ is obtained at any time $t$ in terms of $x^0_{t-1}$ and $z_{t-1}$ using~\eqref{eq:approximation process}.  In addition, for any $i \in \mathbb{N}_n$, follower $i$ solves two Riccati equations~\eqref{eq:Riccati-1} and~\eqref{eq:Riccati-2} to find $\breve L_t, \bar L^{2,1}_t$ and $\bar L^{2,2}_t$, and then computes its control action based on~\eqref{eq:our strategy for followers} using its local state $x^i_t$, the state of the leader $x^0_t$, and variable $z_t$. The results are depicted in Figure~\ref{fig 1}, where the thick curve represents the state of  the leader, and thin curves are the states of the followers (to avoid a cluttered figure, only 100 followers are chosen, randomly, to display their states). It can be observed from this figure that the states of all agents (which are, in fact, the position of the agents) are convergent under the proposed control strategy, as expected, and hence consensus is achieved asymptotically.    

The next example demonstrates the efficacy of the results obtained in this work, for the special case of leaderless multi-agent networks. 

\textbf{Example 2:} Consider a multi-agent system consisting of $100$ agents that are to track a constant reference trajectory $x^0_1=3$.  The following parameters are  used in the simulation:
\begin{align}
& A^0_t=1, \quad B^0_t=0, \quad D^0_t=0,  \quad A_t=1, \quad  B_t=0.5, \\ 
& D_t=0.05, \quad E_t=0,  \quad Q^0_t=0, \quad R^0_t=0,\\  
&Q_t=0.1,  \quad  P_t=20, \quad  R_t=100, \quad H_t=0.5, \quad F_t=60,\\
&T=40, \quad w^i_t \sim \mathcal{N}(0,0.02) \quad \forall i \in \mathbb{N}_{100}.
\end{align}
Similar to Example 1, each follower computes its control action according to \eqref{eq:our strategy for followers}. It is to be noted that  in the leaderless case, the agents do not  communicate as discussed in~Remark~\ref{remark:leaderless}.
\begin{figure}[h]
	\centering
	\vspace{-3.8cm}
	\includegraphics[width=\linewidth]{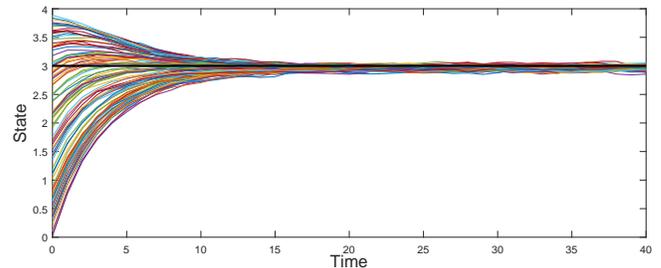} 
	\vspace{-4.2cm}
	\caption{Trajectories of the followers  in Example 2.}\label{fig 2}
\end{figure}
The results are given in Figure~\ref{fig 2}, analogously to Figure~\ref{fig 1}, and show that consensus is achieved as the states of all followers converge to the same value.

\section{Conclusions}\label{sec:conclusion}
A mean-field approach to the decentralized control of a leader-follower multi-agent network with a single leader is presented in this paper, where the states of the leader and followers are coupled in the dynamics and cost. A near-optimal strategy for a non-classical information structure is proposed such that the strategy is obtained by solving two decoupled Riccati equations, where the dimension of the matrices in these equations is independent of the number of followers. This means that the proposed method is not only distributed, it is also scalable. It is shown that the proposed solution converges to the optimal strategy at a rate inversely proportional to the number of followers. The effectiveness of the results is verified by simulation, for two different multi-agent settings with 1000 and 100 followers.

As suggestions for future research directions, one can extend the results to the case of infinite horizon, multiple leaders, heterogeneous followers, and weighted cost functions, under standard assumptions in mean-field teams~\cite{arabneydi2016new}. The approach is robust in the sense that the failure of a small number of followers has negligible impact on the mean-field for a network of sufficiently large population.

\bibliographystyle{IEEEtran}

\bibliography{Jalal_Ref} 
\end{document}